\documentclass[12pt]{elsart}
\usepackage{indentfirst,latexsym}
\usepackage{latexsym}
\usepackage{indentfirst}
\usepackage{flafter}
\usepackage{graphicx}
\usepackage{amssymb}
\newtheorem{theorem}{Theorem}[section]

\newtheorem{example}{Example}[section]
\newtheorem{lemma}{Lemma}[section]

\newtheorem{corollary}{Corollary}[section]

\begin{document}
\runauthor{Qian-Ping Guo, Hou-Biao Li, Jin-Song Leng etc.}
\begin{frontmatter}
\title{Some new bounds for the Hadamard product and the Fan product of matrices\thanksref{b}}
\thanks[b]{\small Supported by National Natural Science Foundation of China (11101071, 11271001, 51175443) and the Fundamental Research
Funds for China Scholarship Council.}

\author[]{Qian-Ping Guo},
\corauth[]{Corresponding author.}\ead{lengjs@uestc.edu.cn, lihoubiao0189@163.com or guoqianpinglei@163.com}
\author[]{Hou-Biao Li},
\author[]{Jin-Song Leng}
\address {School of Mathematical Sciences, University of Electronic Science and Technology of
China, Chengdu, 611731, P. R. China}

\begin{abstract}
If $A$ and $B$ are nonnegative matrices, a sharp upper bound on the spectral radius $\rho(A\circ B)$ for the Hadamard product of two 
nonnegative matrices is given, and the minimum eigenvalue $\tau(A\star B)$ of the Fan product of two $M$-matrices $A$ and $B$ is discussed.
In addition, we also give a sharp lower bound on $\tau(A\circ B^{-1})$ for the Hadamard product of $A$ and $B^{-1}$. Several examples, illustrating that the given bound is stronger than the existing bounds, are also given.
\end{abstract}


\end{frontmatter}

\setcounter{equation}{0}
\renewcommand{\theequation}{1.\arabic{equation}}

\section{Introduction}

In this paper, for a positive integer $n$, $N$ denotes the set $\{1,2,\cdots,n\}$. $\mathbb{R}^{n\times n}$ denotes the set of all $n\times n$ real matrices and the set of all $n\times n$ complex
matrices is denoted by $\mathbb{C}^{n\times n}$. Let $A = (a_{ij})$ and $B = (b_{ij})$ be two real $n\times n$ matrices. We write $A\geq B (> B)$ if $a_{ij} \geq b_{ij} (> b_{ij})$ for all $i,j\in N$. If $A\geq 0(> 0)$, we say that $A$ is a nonnegative (positive) matrix. The spectral radius of $A$ is denoted by $\rho(A)$. If $A$ is a nonnegative matrix, the Perron-Frobenius theorem guarantees that $\rho(A)\in \sigma(A)$,  where $\sigma(A)$ is the set of all eigenvalues of $A$. In addition, define
$\tau(A)\triangleq min\{\lambda|\lambda\in\sigma(A)\}$, and denote by $\mathcal{M}_{n}$ the set of nonsingular $M$-matrices
(see \cite{QG}).

For $n\geqslant 2$, an $n\times n$ matrix $A$ is said to be
reducible if there exists a permutation matrix $P$ such that
$$
\mathop {\mathop {P^{T}AP = \left( {\begin{array}{*{20}c}
   {B} & {C}   \\
   {0} & {D}   \\
\end{array}} \right)}\limits_{} },
$$
where $B$ and $D$ are square matrices of order at least one. If no
such permutation matrix exists, then A is called irreducible. If $A$
is a $1\times1$ complex matrix, then $A$ is irreducible if and only
if its single entry is nonzero (see \cite{GQ}).

According to Ref. \cite{GQ}, a matrix $A$ is called an $M$-matrix,
if there exists an $n\times n$ nonnegative real matrix $P$ and a
nonnegative real number $\alpha$ such that $A=\alpha I-P$, and
$\alpha\geq\rho(P)$, where $\rho(P)$ denotes the spectral radius of
$P$ and $I$ is the identity matrix. Moreover, if $\alpha>\rho(P)$,
$A$ is called a nonsingular $M$-matrix; if $\alpha=\rho(P)$, we call
$A$ a singular $M$-matrix.

In addition, a matrix $A=(a_{ij})\in \mathbb{R}^{n\times n}$ is
called $Z$-matrix if all of whose off-diagonal entries are negative,
and denoted by $A\in \mathcal{Z}_{n}$. For convenience, the following simple
facts are needed (see Problems 16, 19 and 28 in Section 2.5 of
\cite{R.C}):

\begin{enumerate}
  \item $\tau(A)\in\sigma(A)$;
  \item If $A, B\in \mathcal{M}_{n}$, and $A\geq B$, then $\tau(A)\geq\tau(B)$;
  \item If $A\in \mathcal{M}_{n}$, then $\rho(A^{-1})$ is the Perron eigenvalue of the nonnegative matrix $A^{-1}$, and $\tau(A)=\frac{1}{\rho(A^{-1})}$ is a positive real eigenvalue of $A$.
\end{enumerate}

Let $A$ be an irreducible nonsingular $M$-matrix. It is well known
that there exist positive vectors $u$ and $v$ such that
$Au=\tau(A)u$ and $v^{T}A=\tau(A)v^{T}$, where $u$ and $v$ are right
and left Perron eigenvectors of $A$, respectively.

The Hadamard product of $A=(a_{ij})\in \mathbb{C}^{n\times n}$ and $B=(b_{ij})\in \mathbb{C}^{n\times n}$ is defined by
$A\circ B=(a_{ij}b_{ij})\in \mathbb{C}^{n\times n}$.

For two real matrices $A, B\in \mathcal{M}_{n}$, the Fan product of $A$ and $B$ is denoted by $A\star B=C=[c_{ij}]\in \mathcal{M}_{n}$ and is defined by
$$
c_{ij}=\left\{{\begin{array}{ll}
-a_{ij}b_{ij},&{\rm if}\;\;i\neq j,\\
a_{ii}b_{ii},&{\rm if}\;\;i=j.
\end{array}}\right.
$$
We define: for any $i,j,l\in N$,
$$r_{li}=\frac{|a_{li}|}{|a_{ll}|-\sum\limits_{k\neq l,i}|a_{lk}|},~~l\neq i;~~~~r_{i}=\max_{l\neq i}\{r_{li}\},~~i\in N,
$$
$$
s_{ji}=\frac{|a_{ji}|+\sum\limits_{k\neq j,i}|a_{jk}|r_{k}}{|a_{jj}|},~~j\neq i;~~~~s_{i}=\max\limits_{j\neq i}\{s_{ji}\}, ~~i\in N,
$$

throughout the paper.

For two nonnegative matrices $A, B$, we will exhibit a new upper bound for $\rho(A\circ B)$, a new lower bound on the eigenvalue $\tau(A\star B)$ for the Fan product and a new lower bound on the eigenvalue $\tau(A\circ B^{-1})$ for the hadamard product in this paper.

\setcounter{equation}{0}
\renewcommand{\theequation}{2.\arabic{equation}}
 \section{An upper bound for the spectral radius of the Hadamard product of two nonnegtive matrices}\label{sec2}

 In (\cite{R.C}, p. 358), there is a simple estimate for $\rho(A\circ B)$: if $A, B\in \mathbb{R}^{n\times n}$, $A\geq 0$, and $B\geq 0$, then
 \begin{equation}\label{eq:2.1}
 \rho(A\circ B)\leq \rho(A)\rho(B).
\end{equation}
Fang \cite{MF} gave an upper bound for $\rho(A\circ B)$, that is,
\begin{equation}\label{eq:2.2}
\rho(A\circ B)\leq \max_{1\leq i\leq n}\Big\{2a_{ii}b_{ii}+\rho(A)\rho(B)-b_{ii}\rho(A)-a_{ii}\rho(B)\Big\},
\end{equation}
which is shaper than the bound $\rho(A)\rho(B)$ in (\cite{R.C}, p. 358).

Recently, Liu \cite{QG} improved the above results, have
\begin{equation}\label{eq:2.3}
\begin{array}{lll}
\rho(A\circ B)\leq \max\limits_{i\neq j}\frac{1}{2}\Big\{a_{ii}b_{ii}+a_{jj}b_{jj}+[(a_{ii}b_{ii}-a_{jj}b_{jj})^{2}\\
~~~~~~~~~~~~~~+4(\rho(A)-a_{ii})(\rho(B)-b_{ii})(\rho(A)-a_{jj})(\rho(B)-b_{jj})]^{\frac{1}{2}}\Big\}.
\end{array}
\end{equation}

Firstly, we give some lemmas in this section.

\begin{lemma}\label{lem:2.1}(Perron-Frobenius theorem)(\cite{R.C}).
If $A$ is an irreducible nonnegative matrix, there exist positive vectors $u$, such that $Au=\rho(A)u$.
\end{lemma}

\begin{lemma}\label{lem:2.2}(\cite{R.C}).
If $A, B\in \mathbb{C}^{n\times n}$, $D$ and $E$ are positive diagonal matrices, then $$D(A\circ B)E=(DAE)\circ B=(DA)\circ (BE)=(AE)\circ (DB)=A\circ (DBE).$$
\end{lemma}

\begin{lemma}\label{lem:2.3}(Brauer's theorem).
Let $A=(a_{ij})\in\mathbb{C}^{n\times n}$ $(n\geq2)$,
then all the eigenvalues of $A$ lie inside the union of $\frac{n(n-1)}{2}$ ovals of Cassini, i.e.,
\begin{equation}\label{eq:2.4}
B(A)=\bigcup^{n}_{i,j=1;i\neq j}\Big\{z\in
\mathbb{C}:|z-a_{ii}||z-a_{jj}|\leq (\sum\limits_{k\neq
i}|a_{ki}|)(\sum\limits_{k\neq
j}|a_{kj}|)\Big\},
\end{equation}
\end{lemma}

Obviously, if we denote $C=D^{-1}AD$, $D=diag(d_{1},d_{2},\cdots,d_{n}), d_{i}>0$, then $C$ and $A$ have the same eigenvalues, we obtain that all the eigenvalues of $A$ lie in the region:

\begin{equation}\label{eq:2.5}
\bigcup^{n}_{i,j=1;i\neq j}\Big\{z\in
\mathbb{C}:|z-a_{ii}||z-a_{jj}|\leq (\sum\limits_{k\neq
i}\frac{d_{k}}{d_{i}}|a_{ik}|)(\sum\limits_{k\neq
j}\frac{d_{l}}{d_{j}}|a_{jl}|)\Big\}.
\end{equation}

Next, we present a new estimating formula on the upper bound of $\rho(A\circ B)$.

\begin{theorem}\label{th:2.1}
If $A=(a_{ij})$ and $B=(b_{ij})$ are nonnegative matrices, $s_{i}=\max\limits_{j\neq i}\{a_{ij}\}$, $t_{i}=\max\limits_{j\neq i}\{b_{ij}\}$, then
\begin{equation}\label{eq:2.6}
\begin{array}{lll}
\rho(A\circ B)\leq \max \limits_{i\neq j} \frac{1}{2}\Big\{a_{ii}b_{ii}+a_{jj}b_{jj}+[(a_{ii}b_{ii}-a_{jj}b_{jj})^{2}\\
~~~~~~~~~~~~~~~~~~+4t_{i}s_{j}(\rho(A)-a_{ii})(\rho (B)-b_{jj})]^{\frac{1}{2}}\Big\}.
\end{array}
\end{equation}
\end{theorem}

\textbf{Proof.} It is evident that the inequality (\ref{eq:2.6})
holds with the equality for $n=1$. Therefore, we assume that $n\geq 2$
and divide two cases to prove this problem.

\textbf{Case 1}. Suppose that $A\circ B$ is irreducible. Obviously
$A$ and $B$ are also irreducible. By Lemma \ref{lem:2.1}, there
exists positive vectors $u=(u_{1},u_{2},\cdots,u_{n})$ and have $$(D^{-1}AD)u=\rho(D^{-1}AD)u=\rho(A)u,$$ where $D=diag(d_{1},d_{2},\cdots,d_{n}), d_{i}>0$, then $$\sum\limits_{j\neq i}\frac{a_{ij}d_{j}u_{j}}{d_{i}u_{i}}=\rho(A)-a_{ii}.$$

Define $U=diag(u_{1},u_{2},\cdots,u_{n})$, $C=(DU)^{-1}A(DU)$, then we have that
$$
\mathop {\mathop {C=\left( {\begin{array}{*{20}c}
a_{11}&\frac{d_{2}u_{2}}{d_{1}u_{1}}a_{12}&\cdots&\frac{d_{n}u_{n}}{d_{1}u_{1}}a_{1n}\\
\frac{d_{1}u_{1}}{d_{2}u_{2}}a_{21}&a_{22}&\cdots&\frac{d_{n}u_{n}}{d_{2}u_{2}}a_{2n}\\
\vdots&\vdots&\ddots&\vdots\\
\frac{d_{1}u_{1}}{d_{n}u_{n}}a_{n1}&\frac{d_{2}u_{2}}{d_{n}u_{n}}a_{n2}&\cdots&a_{nn}\\
\end{array}} \right)}\limits_{} }
$$
is an irreducible nonnegative matrix and

$$
\mathop {\mathop {C\circ B=(m_{ij})=\left( {\begin{array}{*{20}c}
a_{11}b_{11}&\frac{d_{2}u_{2}}{d_{1}u_{1}}a_{12}b_{12}&\cdots&\frac{d_{n}u_{n}}{d_{1}u_{1}}a_{1n}b_{1n}\\
\frac{d_{1}u_{1}}{d_{2}u_{2}}a_{21}b_{21}&a_{22}b_{22}&\cdots&\frac{d_{n}u_{n}}{d_{2}u_{2}}a_{2n}b_{2n}\\
\vdots&\vdots&\ddots&\vdots\\
\frac{d_{1}u_{1}}{d_{n}u_{n}}a_{n1}b_{n1}&\frac{d_{2}u_{2}}{d_{n}u_{n}}a_{n2}b_{n2}&\cdots&a_{nn}b_{nn}\\
\end{array}} \right)}\limits_{} }.
$$

By Lemma \ref{lem:2.2},
$$(DU)^{-1}(A\circ B)(DU)=(DU)^{-1}A(DU)\circ B=C\circ B,$$ i.e., $\rho (A\circ B)=\rho (C\circ B)$.

By the inequality (\ref{eq:2.4}) and $\rho(A\circ B)\geq a_{ii}b_{ii}$ (see \cite{A.R}), $\mathrm{for ~any} ~j\neq i\in N$,

 we have
\begin{equation}
\begin{array}{lll}\label{eq:2.7}
(\rho(A\circ B)-a_{ii}b_{ii})(\rho(A\circ B))-a_{jj}b_{jj})
\leq \sum\limits_{ k\neq i}|m_{ik}|\sum\limits_{l\neq j}|m_{jl}|\\
~~~~~~~~~~~~~~~~~~~~~~~=\sum\limits_{k\neq i}\frac{d_{k}u_{k}a_{ik}b_{ik}}{d_{i}u_{i}}\sum\limits_{l\neq j}\frac{d_{l}u_{l}a_{jl}b_{jl}}{d_{j}u_{j}}\\
~~~~~~~~~~~~~~~~~~~~~~~\leq \Big(\max\limits_{k\neq i}\{b_{ik}\}\sum\limits_{k\neq i}\frac{d_{k}u_{k}a_{ik}}{d_{i}u_{i}}\Big)\Big(\max\limits_{l\neq j}\{a_{jl}\}\sum\limits_{l\neq j}\frac{d_{l}u_{l}b_{jl}}{d_{j}u_{j}}\Big)\\
~~~~~~~~~~~~~~~~~~~~~~~\leq\max\limits_{k\neq i}\{b_{ik}\}(\rho (A)-a_{ii})\max\limits_{l\neq j}\{a_{jl}\}(\rho (B)-b_{jj})\\
~~~~~~~~~~~~~~~~~~~~~~~=t_{i}s_{j}\Big(\rho (A)-a_{ii}\Big)\Big(\rho (B)-b_{jj}\Big).\\
\end{array}
\end{equation}
Thus, by solving the quadratic inequality (\ref{eq:2.7}), we have that
$$
\begin{array}{lll}
\rho(A\circ B)\leq \frac{1}{2} \Big\{a_{ii}b_{ii}+a_{jj}b_{jj}+[(a_{ii}b_{ii}-a_{jj}b_{jj})^{2}
+4t_{i}s_{j}(\rho (A)-a_{ii})(\rho (B)-b_{jj})]^{\frac{1}{2}}\Big\}\\
~~~~~~~~~~~\leq\max\limits_{i\neq
j}\frac{1}{2} \Big\{a_{ii}b_{ii}+a_{jj}b_{jj}+[(a_{ii}b_{ii}-a_{jj}b_{jj})^{2}
+4t_{i}s_{j}(\rho (A)-a_{ii})(\rho (B)-b_{jj})]^{\frac{1}{2}}\Big\}.
\end{array}
$$
i.e., the conclusion (\ref{eq:2.6}) holds.

\textbf{Case 2}. If $A\circ B$ is reducible. We may denote by $P=(p_{ij})$ the
$n\times n$ permutation matrix $(p_{ij})$ with
$$p_{12}=p_{23}=\cdots=p_{n-1,n}=p_{n,1}=1,$$ the remaining $p_{ij}$ zero, then both $A+\varepsilon P$ and $B+\varepsilon P$ are nonnegative irreducible matrices for any sufficiently small positive real number $\varepsilon$. Now we substitute $A+\varepsilon P$ and $B+\varepsilon P$
for $A$ and $B$, respectively in the previous Case 1, and then
letting $\varepsilon\rightarrow 0$, the result (\ref{eq:2.6})
follows by continuity. $\Box$

\textbf{Remark 2.1.}\label{re:2.1} Next, we give a comparison between
the upper bound in the inequality (\ref{eq:2.3}) and the upper bound in
the inequality (\ref{eq:2.6}). Without loss of generality, if $t_{i}+b_{ii}\geq\rho(B)$, $s_{j}+a_{jj}\geq\rho(A)$, $i,j=1,\cdots,n$, then we have $t_{i}s_{j}\geq(\rho(B)-b_{ii})(\rho(A)-a_{jj})$.
Thus, the upper bound in the inequality (\ref{eq:2.6}) is better than the upper bound in the inequality (\ref{eq:2.3}).

\begin{example}\label{ex:2.1}.
Let A and B be the same as in Example 1 from \cite{QG}:
\end{example}

$$
\mathop {\mathop {A=(a_{ij})=\left( {\begin{array}{*{20}c}
{4}&{1}&{0}&{2}\\
{1}&{0.05}&{1}&{1}\\
{0}&{1}&{4}&{0.5}\\
{1}&{0.5}&{0}&{4}
\end{array}} \right)}\limits_{} }, ~~
\mathop {\mathop {B=(b_{ij})=\left( {\begin{array}{*{20}c}
{1}&{1}&{1}&{1}\\
{1}&{1}&{1}&{1}\\
{1}&{1}&{1}&{1}\\
{1}&{1}&{1}&{1}
\end{array}} \right).}\limits_{} }
$$

By direct calculation, $\rho(A\circ B)=5.7339$.

According to (\ref{eq:2.1}), we have

$$\rho(A\circ B)\leq \rho(A)\rho(B)=22.9336.$$

If we apply (\ref{eq:2.2}) and (\ref{eq:2.3}), we get

$$\rho(A\circ B)\leq \max\limits_{1\leq i\leq 4}\Big\{2a_{ii}b_{ii}+\rho(A)\rho(B)-a_{ii}\rho(B)-b_{ii}\rho(A)\Big\}=17.1017,$$

and

$$
\begin{array}{lll}
\rho(A\circ B)\leq \max\limits_{i\neq j}\frac{1}{2}\Big\{a_{ii}b_{ii}+a_{jj}b_{jj}+[(a_{ii}b_{ii}-a_{jj}b_{jj})^{2}\\
~~~~~~~~~~~~~~+4(\rho(A)-a_{ii})(\rho(B)-b_{ii})(\rho(A)-a_{jj})(\rho(B)-b_{jj})]^{\frac{1}{2}}\Big\}=11.6478.
\end{array}
$$

If we apply Theorem \ref{th:2.1}, we obtain that

$$
\begin{array}{lll}
\rho(A\circ B)\leq \max \limits_{i\neq j} \frac{1}{2}\Big\{a_{ii}b_{ii}+a_{jj}b_{jj}+[(a_{ii}b_{ii}-a_{jj}b_{jj})^{2}\\
~~~~~~~~~~~~~~~~~~+4t_{i}s_{j}(\rho(A)-a_{ii})(\rho (B)-b_{jj})]^{\frac{1}{2}}\Big\}=8.1897.
\end{array}
$$

The example shows that the bound in Theorem \ref{th:2.1} is better than the existing bounds.

In addition, by the Theorem \ref{th:2.1} and \cite{QG}, we also have the following corollary:

 \begin{corollary}\label{cor:2.1}
Let A and B be nonnegative matrices, then
\end{corollary}

 $$
\begin{array}{lll}
|det(A\circ B)|\leq \Big(\rho(A\circ B)\Big)^{n}\\
~~~~~~~~~~~~~~~~\leq\max\limits_{i\neq j} \frac{1}{2^{n}}\Big\{a_{ii}b_{ii}+a_{jj}b_{jj}+[(a_{ii}b_{ii}-a_{jj}b_{jj})^{2}\\
~~~~~~~~~~~~~~~~~~~~~~+4t_{i}s_{j}(\rho(A)-a_{ii})(\rho (B)-b_{jj})]^{\frac{1}{2}}\Big\}^{n} \\
~~~~~~~~~~~~~~~~~\leq\max \limits_{i\neq j} \frac{1}{2^{n}}\Big\{a_{ii}b_{ii}+a_{jj}b_{jj}+[(a_{ii}b_{ii}-a_{jj}b_{jj})^{2}\\
~~~~~~~~~~~~~~~~~~~~~~+4(\rho(A)-a_{ii})(\rho (B)-b_{ii})(\rho(A)-a_{jj})(\rho (B)-b_{jj})]^{\frac{1}{2}}\Big\}^{n}.
\end{array}
$$

\setcounter{equation}{0}
\renewcommand{\theequation}{3.\arabic{equation}}
\section{Inequalities for the Fan product of two $M$-matrices}\label{sec3}

It is known (p.359, \cite{R.C}) that the following classical result
is given: if $A,B\in \mathbb{R}^{n\times n}$ are $M$-matrices, then
\begin{equation}\label{eq:3.1}
 \tau(A\star B)\geq\tau(A)\tau(B).
\end{equation}
In 2007, Fang improved (\ref{eq:3.1}) in the Remark 3 of Ref.
\cite{MF} and gave a new lower bound for
 $\tau(A\star B)$, that is
 \begin{equation}\label{eq:3.2}
 \tau(A\star B)\geq \min_{1\leq i\leq n}\Big\{b_{ii}\tau(A)+a_{ii}\tau(B)-\tau(A)\tau(B)\Big\}.
 \end{equation}

Subsequently, Liu et al.\cite{QG} gave a sharper bound than
(\ref{eq:3.2}), i.e.,
 \begin{equation}
\begin{array}{lll}\label{eq:3.3}
\tau(A\star B)\geq \frac{1}{2} \min \limits_{i\ne j}\Big\{a_{ii}b_{ii}+a_{jj}b_{jj}-[(a_{ii}b_{ii}-a_{jj}b_{jj})^{2}\\
~~~~~~~~~~~~~~+4(b_{ii}-\tau(B))(a_{ii}-\tau(A))(b_{jj}-\tau(B))(a_{jj}-\tau(A))]^{\frac{1}{2}}\Big\}.
\end{array}
\end{equation}

In addition, by the definition of Fan product, the following lemma holds:

\begin{lemma}\label{lem:3.1}(\cite{QG}).
If $A, B\in \mathbb{C}^{n\times n}$ be nonsingular $M$-matrices, $D$ and $E$ are positive diagonal matrices, then $$D(A\star B)E=(DAE)\star B=(DA)\star (BE)=(AE)\star (DB)=A\star (DBE).$$
\end{lemma}
Next, we give a new lower bound on the minimum eigenvalue $\tau(A\star B)$ of the Fan product of nonsingular $M$-matrices.
\begin{theorem}\label{th:3.1}
If $A=(a_{ij})$ and $B=(b_{ij})$ are nonsingular $M$-matrices, $s_{i}=\max\limits_{j\neq i}|a_{ij}|$, $t_{i}=\max\limits_{j\neq i}|b_{ij}|$, then
\begin{equation}\label{eq:3.4}
\begin{array}{lll}
\tau(A\star B)\geq \min \limits_{i\neq j}\frac{1}{2}\Big\{a_{ii}b_{ii}+a_{jj}b_{jj}-[(a_{ii}b_{ii}-a_{jj}b_{jj})^{2}\\
~~~~~~~~~~~~~~~~~~+4t_{i}s_{j}(a_{ii}-\tau (A))(b_{jj}-\tau (B))]^{\frac{1}{2}}\Big\}.
\end{array}
\end{equation}
\end{theorem}

\textbf{Proof.} It is clear that the (\ref{eq:3.4})
holds with the equality for $n=1$.

We next assume $n\geq 2$ and divide two cases to prove this problem.

\textbf{Case 1}. Suppose that $A\star B$ is irreducible. Obviously
$A$ and $B$ are also irreducible. By \cite{A.R}, there
exists positive vectors $u=(u_{1},u_{2},\cdots,u_{n})$ such that $$(D^{-1}AD)u=\tau(D^{-1}AD)u=\tau(A)u,$$ where $D=diag(d_{1},d_{2},\cdots,d_{n}), d_{i}>0$, 
and then $$a_{ii}-\sum\limits_{j\neq i}\frac{|a_{ij}|d_{j}u_{j}}{d_{i}u_{i}}=\tau(A).$$

Define $U=diag(u_{1},u_{2},\cdots,u_{n})$, $C=(DU)^{-1}A(DU)$, we have that

$$
\mathop {\mathop {C=\left( {\begin{array}{*{20}c}
a_{11}&\frac{d_{2}u_{2}}{d_{1}u_{1}}a_{12}&\cdots&\frac{d_{n}u_{n}}{d_{1}u_{1}}a_{1n}\\
\frac{d_{1}u_{1}}{d_{2}u_{2}}a_{21}&a_{22}&\cdots&\frac{d_{n}u_{n}}{d_{2}u_{2}}a_{2n}\\
\vdots&\vdots&\ddots&\vdots\\
\frac{d_{1}u_{1}}{d_{n}u_{n}}a_{n1}&\frac{d_{2}u_{2}}{d_{n}u_{n}}a_{n2}&\cdots&a_{nn}\\
\end{array}} \right)}\limits_{} }
$$

is an irreducible nonsingular $M$- matrix, then

$$
\mathop {\mathop {C\star B=(m_{ij})=\left( {\begin{array}{*{20}c}
a_{11}b_{11}&\frac{d_{2}u_{2}}{d_{1}u_{1}}a_{12}b_{12}&\cdots&\frac{d_{n}u_{n}}{d_{1}u_{1}}a_{1n}b_{1n}\\
\frac{d_{1}u_{1}}{d_{2}u_{2}}a_{21}b_{21}&a_{22}b_{22}&\cdots&\frac{d_{n}u_{n}}{d_{2}u_{2}}a_{2n}b_{2n}\\
\vdots&\vdots&\ddots&\vdots\\
\frac{d_{1}u_{1}}{d_{n}u_{n}}a_{n1}b_{n1}&\frac{d_{2}u_{2}}{d_{n}u_{n}}a_{n2}b_{n2}&\cdots&a_{nn}b_{nn}\\
\end{array}} \right)}\limits_{} }.
$$

By the Lemma \ref{lem:3.1},
$$(DU)^{-1}(A\star B)(DU)=(DU)^{-1}A(DU)\star B=C\star B,$$ i.e., $\tau (A\star B)=\tau (C\star B)$.

In addition, by the inequality (\ref{eq:2.4}) and $0\leq\tau(A\star B)\leq a_{ii}b_{ii}$ (see \cite{A.R}), $\mathrm{for ~any} ~j\neq i\in N$, we have
\begin{equation}
\begin{array}{lll}\label{eq:3.5}
|\tau(A\star B)-a_{ii}b_{ii}||\tau(A\star B)-a_{jj}b_{jj}|
\leq \sum\limits_{ k\neq i}|m_{ik}|\sum\limits_{l\neq j}|m_{jl}|\\
~~~~~~~~~~~~~~~~~~~~~~~=\sum\limits_{k\neq i}|\frac{d_{k}u_{k}a_{ik}b_{ik}}{d_{i}u_{i}}|\sum\limits_{l\neq j}|\frac{d_{l}u_{l}a_{jl}b_{jl}}{d_{j}u_{j}}|\\
~~~~~~~~~~~~~~~~~~~~~~~\leq \Big(\max\limits_{k\neq i}|b_{ik}|\sum\limits_{k\neq i}|\frac{d_{k}u_{k}a_{ik}}{d_{i}u_{i}}|\Big)\Big(\max\limits_{l\neq j}|a_{jl}|\sum\limits_{l\neq j}|\frac{d_{l}u_{l}b_{jl}}{d_{j}u_{j}}|\Big)\\
~~~~~~~~~~~~~~~~~~~~~~~\leq\max\limits_{k\neq i}|b_{ik}|(a_{ii}-\tau (A))\max\limits_{l\neq j}|a_{jl}|(b_{jj}-\tau (B))\\
~~~~~~~~~~~~~~~~~~~~~~~=t_{i}s_{j}(a_{ii}-\tau (A))(b_{jj}-\tau (B)).\\
\end{array}
\end{equation}
Thus, by solving the quadratic inequality (\ref{eq:3.5}), we have that
$$
\begin{array}{lll}
\tau(A\star B)\geq \frac{1}{2} \Big\{a_{ii}b_{ii}+a_{jj}b_{jj}-[(a_{ii}b_{ii}-a_{jj}b_{jj})^{2}
+4t_{i}s_{j}(a_{ii}-\tau (A))(b_{jj}-\tau (B))]^{\frac{1}{2}}\Big\}\\
~~~~~~~~~~~\geq\min\limits_{i\neq
j}\frac{1}{2} \Big\{a_{ii}b_{ii}+a_{jj}b_{jj}-[(a_{ii}b_{ii}-a_{jj}b_{jj})^{2}
+4t_{i}s_{j}(a_{ii}-\tau (A))(b_{jj}-\tau (B))]^{\frac{1}{2}}\Big\}.
\end{array}
$$
i.e., the conclusion (\ref{eq:3.4}) holds.

\textbf{Case 2}. If $A\star B$ is reducible. It is well known that a matrix in $\mathcal{Z}_{n}$ is a nonsingular $M$-matrix if and only if all its leading principal minors are positive (see condition (E17) of Theorem 6.2.3 of \cite{A.R}). We denote by $P=(p_{ij})$ the $n\times n$ permutation matrix $(p_{ij})$ with
$$p_{12}=p_{23}=\cdots=p_{n-1,n}=p_{n,1}=1,$$ the remaining $p_{ij}$ zero, then both $A-\varepsilon P$ and $B-\varepsilon P$ are irreducible nonsingular $M$-matrices for any sufficiently small positive real number $\varepsilon$. Now we substitute $A-\varepsilon P$ and $B-\varepsilon P$ for $A$ and $B$, respectively in the previous Case 1, and then letting $\varepsilon\rightarrow 0$, the result (\ref{eq:3.4}) follows by continuity. $\Box$

\textbf{Remark 3.1.}\label{re:3.1} Similarly, we give a comparison between
the lower bound in the inequality (\ref{eq:3.3}) and the lower bound in
the inequality (\ref{eq:3.4}).
If $a_{jj}\geq\tau(A)+s_{j}$, $b_{ii}\geq\tau(B)+t_{i}$, $i,j=1,\cdots,n$, then $(a_{jj}-\tau(A))(b_{ii}-\tau(B))\geq s_{j}t_{i}$ for all $i\neq j$.
Thus, the lower bound in the inequality (\ref{eq:3.4}) is better than the lower bound in the inequality (\ref{eq:3.3}).

In addition, from Theorem \ref{th:3.1} and \cite{A.R}, we may get the following corollary.

\begin{corollary}\label{cor:3.1}.
If $A$, $B$ are nonsingular $M$-matrices, then
\end{corollary}
$$
\begin{array}{lll}
|det(A\star B)|\geq \Big(\tau(A\star B)\Big)^{n}\\
~~~~~~~~~~~~~~~~\geq\min \limits_{i\neq j}\frac{1}{2^{n}}\Big\{a_{ii}b_{ii}+a_{jj}b_{jj}-[(a_{ii}b_{ii}-a_{jj}b_{jj})^{2}\\
~~~~~~~~~~~~~~~~~~+4t_{i}s_{j}(a_{ii}-\tau (A))(b_{jj}-\tau (B))]^{\frac{1}{2}}\Big\}^{n} \\
~~~~~~~~~~~~~~~~~\geq\min \limits_{i\neq j}\frac{1}{2^{n}}\Big\{a_{ii}b_{ii}+a_{jj}b_{jj}-[(a_{ii}b_{ii}-a_{jj}b_{jj})^{2}\\
~~~~~~~~~~~~~~~~~~+4(a_{ii}-\tau (A))(b_{ii}-\tau (B))(a_{jj}-\tau (A))(b_{jj}-\tau (B))]^{\frac{1}{2}}\Big\}^{n}.
\end{array}
$$

\begin{example}\label{ex:3.1}(\cite{QG}).
Let A and B be the nonsingular $M$-matrices:
\end{example}

$$
\mathop {\mathop {A=(a_{ij})=\left( {\begin{array}{*{20}c}
{2}&{-1}&{0}\\
{0}&{1}&{-0.5}\\
{-0.5}&{-1}&{2}
\end{array}} \right)}\limits_{} }, ~~
\mathop {\mathop {B=(b_{ij})=\left( {\begin{array}{*{20}c}
{1}&{-0.25}&{-0.25}\\
{-0.5}&{1}&{-0.25}\\
{-0.25}&{-0.5}&{1}
\end{array}} \right).}\limits_{} }
$$

By (\ref{eq:3.1}), we have

$$\tau(A\star B)\geq \tau(A)\tau(B)=0.1854.$$

If we use the inequalities (\ref{eq:3.2}) and (\ref{eq:3.3}), then we get

$$\tau(A\star B)\geq \min\limits_{1\leq i\leq 3}\Big\{a_{ii}\tau(B)+b_{ii}\tau(A)-\tau(A)\tau(B)\Big\}=0.6980,$$

and
$$
\begin{array}{lll}
\tau(A\star B)\geq \min\limits_{i\neq j}\frac{1}{2}\Big\{a_{ii}b_{ii}+a_{jj}b_{jj}-[(a_{ii}b_{ii}-a_{jj}b_{jj})^{2}\\
~~~~~~~~~~~~~~+4(a_{ii}-\tau(A))(b_{ii}-\tau(B))(a_{jj}-\tau(A))(b_{jj}-\tau(B))]^{\frac{1}{2}}\Big\}=0.7655.
\end{array}
$$

If we apply Theorem \ref{th:3.1}, we obtain that

$$
\begin{array}{lll}
\tau(A\star B)\geq \min \limits_{i\neq j}\frac{1}{2}\Big\{a_{ii}b_{ii}+a_{jj}b_{jj}-[(a_{ii}b_{ii}-a_{jj}b_{jj})^{2}\\
~~~~~~~~~~~~~~~+4t_{i}s_{j}(a_{ii}-\tau(A))(b_{jj}-\tau(B))]^{\frac{1}{2}}\Big\}=0.8002.
\end{array}
$$

In fact, $\tau(A\star B)=0.8819$. The example shows that the bound in Theorem \ref{th:3.1} is better than the existing bounds.

\setcounter{equation}{0}
\renewcommand{\theequation}{4.\arabic{equation}}
 \section{A bound for the Hadamard product of $M$-matrix and an inverse $M$-matrix}\label{sec4}

Now, we consider the lower bound of $\tau(A\circ B^{-1})$, for $A=(a_{ij}),B=(b_{ij})\in \mathcal{M}_{n}$ and $B^{-1}=(\beta_{ij})$.

Firstly, in \cite{R.C}, Horn and Johnson gave the classical results
\begin{equation}\label{eq:4.1}
\tau(A\circ B^{-1})\geq\tau(A)\min\limits_{1\leq i\leq n}\beta_{ii}.
\end{equation}

Subsequently, Huang \cite{HR} gave new bound for $\tau(A\circ B^{-1})$, that is,
\begin{equation}\label{eq:4.2}
\tau(A\circ B^{-1})\geq \frac{1-\rho(J_{A})\rho(J_{B})}{1+\rho^{2}(J_{B})}\min\limits_{1\leq i\leq n}\frac{a_{ii}}{b_{ii}},
\end{equation}
where $\rho(J_{A})$ and $\rho(J_{B})$ are the spectral radius of the Jacobi iterative matrices $J_{A}$ and $J_{B}$, respectively.

In 2008, Li \cite{LYT} improved the above results as follows.
\begin{equation}\label{eq:4.3}
\tau(A\circ B^{-1})\geq\min\limits_{i}\frac{b_{ii}-s_{i}\sum\limits_{j\neq i}|b_{ji}|}{a_{ii}}.
\end{equation}

Recently, Chen \cite{CFB} improved the result and gave a new lower bound for $\tau(A\circ B^{-1})$:
\begin{equation}\label{eq:4.4}
\begin{array}{lll}
\tau(A\circ B^{-1})\geq \min \limits_{i\neq j} \frac{1}{2}\Big\{a_{ii}\beta_{ii}+a_{jj}\beta_{jj}-[(a_{ii}\beta_{ii}-a_{jj}\beta_{jj})^{2}\\
~~~~~~~~~~~~~~~~~~+4a_{ii}a_{jj}\beta_{ii}\beta_{jj}\rho^{2}(J_{A})\rho^{2}(J_{B})]^{\frac{1}{2}}\Big\}.
\end{array}
\end{equation}

In this section, we give a lower bound of $\tau(A\circ B^{-1})$ for $M$-matrix and inverse $M$-matrix, which improves the above bounds.

\begin{lemma}\label{lem:4.1}(\cite{CJL}).
If $A=(a_{ij})\in\mathcal{M}_{n}$, there exists a positive diagonal matrix $D$ such that $D^{-1}AD$ is a strictly row diagonally dominant matrix.
\end{lemma}

\begin{lemma}\label{lem:4.2}(\cite{CJL}).
If $A=(a_{ij})\in\mathcal{M}_{n}$, and $D=diag(d_{1},d_{2},\cdots,d_{n})$, $d_{i}>0 ~~(i\in N)$, then $D^{-1}AD$ is also an $M$-matrix.
\end{lemma}

\begin{lemma}\label{lem:4.3}(\cite{CJL}).
If $A, B\in\mathcal{M}_{n}$, then $B\circ A^{-1}$ is also an $M$-matrix.
\end{lemma}

\begin{lemma}\label{lem:4.4}(\cite{LYT}).
If $A=(a_{ij})$ be a strictly diagonally dominant $M$-matrix by rows, then for $A^{-1}=(\alpha_{ij})$, we have $$\alpha_{ji}\leq\frac{|a_{ji}|+\sum\limits_{k\neq j,i}|a_{jk}|r_{k}}{a_{jj}}\alpha_{ii}, ~~~~\mathrm{for ~all} ~~j\neq i.$$
\end{lemma}

\begin{theorem}\label{th:4.1}
If $A=(a_{ij})$ and $B=(b_{ij})$ are two nonsingular $M$-matrices and $B^{-1}=(\beta_{ij})$, $s_{i}=\max\limits_{j\neq i}|a_{ij}|$, then
\begin{equation}\label{eq:4.5}
\begin{array}{lll}
\tau(A\circ B^{-1})\geq \min \limits_{i\neq j} \frac{1}{2}\Big\{a_{ii}\beta_{ii}+a_{jj}\beta_{jj}-[(a_{ii}\beta_{ii}-a_{jj}\beta_{jj})^{2}\\
~~~~~~~~~~~~~~~~~~+4s_{i}s_{j}\beta_{ii}\beta_{jj}(a_{ii}-\tau(A))(b_{jj}-\tau (B))]^{\frac{1}{2}}\Big\}.
\end{array}
\end{equation}
\end{theorem}

\textbf{Proof.} If $A$ is an $M$-matrix, by Lemmas (\ref{lem:4.1}-\ref{lem:4.2}), there exists a positive diagonal matrix $D$ such that
$D^{-1}AD$ is a strictly diagonally dominant $M$-matrix by rows.

\textbf{Case 1}. Suppose that $A\circ B^{-1}$ is irreducible. Obviously
$A$ and $B$ are also irreducible. Since $A-\tau(A)I$ is an irreducible nonsingular $M$-matrix, then $a_{ii}-\tau(A)>0$, $\forall i\in N$, and there exists a positive vector
 $u=(u_{1},u_{2},\cdots,u_{n})$ such that $$Au=\tau(A)u, $$ where $u=diag(u_{1},u_{2},\cdots,u_{n}), u_{i}>0$, and then $$a_{ii}+\sum\limits_{j\neq i}\frac{a_{ji}u_{j}}{u_{i}}=\tau(A).$$

Define $U=diag(u_{1},u_{2},\cdots,u_{n})$, $C=U^{-1}AU$, then we have that
$$
\mathop {\mathop {C=(\tilde{a}_{ij})=U^{-1}AU=\left( {\begin{array}{*{20}c}
a_{11}&\frac{a_{12}u_{1}}{u_{2}}&\cdots&\frac{a_{1n}u_{1}}{u_{n}}\\
\frac{a_{21}u_{2}}{u_{1}}&a_{22}&\cdots&\frac{a_{2n}u_{2}}{u_{n}}\\
\vdots&\vdots&\ddots&\vdots\\
\frac{a_{n1}u_{n}}{u_{1}}&\frac{a_{n2}u_{n}}{u_{2}}&\cdots&a_{nn}\\
\end{array}} \right)}\limits_{} }
$$
is an irreducible nonsingular $M$-matrix.

By Lemma \ref{lem:2.2},
$$U^{-1}(A\circ B^{-1})U=(U^{-1}AU)\circ B^{-1}=C\circ B^{-1},$$ i.e., $\tau (A\circ B^{-1})=\tau (C\circ B^{-1})$.

By the inequality (\ref{eq:2.4}) and $0\leq\tau(A\star B)\leq a_{ii}b_{ii}$ (see \cite{A.R}), $\mathrm{for ~any} ~j\neq i\in N$,

we have
\begin{equation}
\begin{array}{lll}\label{eq:4.6}
|\tau(A\circ B^{-1})-a_{ii}\beta_{ii}||\tau(A\circ B^{-1})-a_{jj}\beta_{jj}|
\leq \sum\limits_{ k\neq i}|\tilde{a}_{ki}|\beta_{ki}\sum\limits_{l\neq j}|\tilde{a}_{lj}|\beta_{lj}\\
~~~~~~~~~~~~~~~~~~~~~~~\leq\sum\limits_{k\neq i}|\tilde{a}_{ki}|\frac{b_{ki}+\sum\limits_{u\neq k,i}|b_{ku}|r_{u}}{b_{kk}}\beta_{ii}\sum\limits_{l\neq j}|\tilde{a}_{lj}|\frac{b_{lj}+\sum\limits_{v\neq l,j}|b_{lv}|r_{v}}{b_{ll}}\beta_{jj}\\
~~~~~~~~~~~~~~~~~~~~~~~=\sum\limits_{k\neq i}|\tilde{a}_{ki}|s_{ki}\beta_{ii}\sum\limits_{l\neq j}|\tilde{a}_{lj}|s_{lj}\beta_{jj}\\
~~~~~~~~~~~~~~~~~~~~~~~\leq \sum\limits_{k\neq i}|\tilde{a}_{ki}|s_{i}\beta_{ii}\sum\limits_{l\neq j}|\tilde{a}_{lj}|s_{j}\beta_{jj}\\
~~~~~~~~~~~~~~~~~~~~~~~=\sum\limits_{k\neq i}\frac{|a_{ki}|u_{k}}{u_{j}}s_{i}\beta_{ii}\sum\limits_{l\neq j}\frac{|a_{jl}|u_{j}}{u_{l}}s_{j}\beta_{jj}\\
~~~~~~~~~~~~~~~~~~~~~~~=s_{i}s_{j}\beta_{ii}\beta_{jj}(a_{ii}-\tau (A))(a_{jj}-\tau (A)).\\
\end{array}
\end{equation}
Thus, by solving the quadratic inequality (\ref{eq:4.6}), we obtain that
$$
\begin{array}{lll}
\tau(A\circ B^{-1})\geq \frac{1}{2} \Big\{a_{ii}\beta_{ii}+a_{jj}\beta_{jj}-[(a_{ii}\beta_{ii}-a_{jj}\beta_{jj})^{2}
+4s_{i}s_{j}\beta_{ii}\beta_{jj}(a_{ii}-\tau (A))(a_{jj}-\tau (A))]^{\frac{1}{2}}\Big\}\\
~~~~~~~~~~~\geq\min\limits_{i\neq
j}\frac{1}{2} \Big\{a_{ii}\beta_{ii}+a_{jj}\beta_{jj}-[(a_{ii}\beta_{ii}-a_{jj}\beta_{jj})^{2}
+4s_{i}s_{j}\beta_{ii}\beta_{jj}(a_{ii}-\tau (A))(a_{jj}-\tau (A))]^{\frac{1}{2}}\Big\}.
\end{array}
$$
i.e., the conclusion (\ref{eq:4.5}) holds.

\textbf{Case 2}. If $A\circ B^{-1}$ is reducible, then one denotes by $P=(p_{ij})$ the
$n\times n$ permutation matrix with $$p_{12}=p_{23}=\cdots=p_{n-1,n}=p_{n,1}=1,$$ the remaining $p_{ij}$ zero, then both $A-\varepsilon P$ and $B-\varepsilon P$ are irreducible nonsingular $M$-matrices for any sufficiently small positive real number $\varepsilon$. Now we substitute $A-\varepsilon P$ and $B-\varepsilon P$
for $A$ and $B$, respectively from the previous Case, and then letting $\varepsilon\rightarrow 0$, the result (\ref{eq:2.6}) follows by continuity. $\Box$

\begin{example}\label{ex:4.1}(\cite{CFB}).
Let A and B be nonsingular $M$-matrices:
\end{example}

$$
\mathop {\mathop {A=(a_{ij})=\left( {\begin{array}{*{20}c}
{1}&{-0.5}&{0}&{0}\\
{-0.5}&{1}&{-0.5}&{0}\\
{0}&{-0.5}&{1}&{-0.5}\\
{0}&{0}&{-0.5}&{1}
\end{array}} \right)}\limits_{} }, ~~
\mathop {\mathop {B=(b_{ij})=\left( {\begin{array}{*{20}c}
{4}&{-1}&{-1}&{-1}\\
{-2}&{5}&{-1}&{-1}\\
{0}&{-2}&{4}&{-1}\\
{-1}&{-1}&{-1}&{4}
\end{array}} \right).}\limits_{} }
$$

By direct calculation, $\tau(A\circ B^{-1})=0.2148$.

According to (\ref{eq:4.1}), we have

$$\tau(A\circ B^{-1})\geq\tau(A)\min\limits_{1\leq i\leq n}\beta_{ii}=0.07.$$

If we apply (\ref{eq:4.2}) and (\ref{eq:4.3}), we get
$$\tau(A\circ B^{-1})\geq\frac{1-\rho(J_{A}\rho(J_{B}))}{1+\rho^{2}(J_{B})}\min\limits_{i}\frac{b_{ii}}{a_{ii}}=0.0707,$$
and
$$\tau(A\circ B^{-1})\geq\min\limits_{i}\frac{b_{ii}-s_{i}\sum\limits_{j\neq i}|b_{ji}|}{a_{ii}}=0.08.$$

According to (\ref{eq:4.4})
$$
\begin{array}{lll}
\tau(A\circ B^{-1})\geq \min \limits_{i\neq j} \frac{1}{2}\Big\{a_{ii}\beta_{ii}+a_{jj}\beta_{jj}-[(a_{ii}\beta_{ii}-a_{jj}\beta_{jj})^{2}\\
~~~~~~~~~~~~~~~~~~+4a_{ii}a_{jj}\beta_{ii}\beta_{jj}\rho^{2}(J_{A})\rho^{2}(J_{B})]^{\frac{1}{2}}\Big\}=0.1524.
\end{array}
$$

If we apply Theorem \ref{th:4.1}, we obtain that

$$
\begin{array}{lll}
\tau(A\circ B^{-1})\geq\min\limits_{i\neq
j}\frac{1}{2} \Big\{a_{ii}\beta_{ii}+a_{jj}\beta_{jj}+[(a_{ii}\beta_{ii}-a_{jj}\beta_{jj})^{2}\\
~~~~~~~~~~~~~~~~+4s_{i}s_{j}\beta_{ii}\beta_{jj}(a_{ii}-\tau (A))(a_{jj}-\tau (A))]^{\frac{1}{2}}\Big\}=0.1929.
\end{array}
$$

 The example shows that the bound in Theorem \ref{th:4.1} is better than the existing bounds.

\setcounter{equation}{0}
\renewcommand{\theequation}{5.\arabic{equation}}
 \section{Inequalities for the Fan product of several $M$-matrices}\label{sec5}

Firstly, let us recall the following lemmas.

\begin{lemma}\label{lem:5.1}(\cite{LHB}).
Let $A$ be an irreducible nonsingular $M$-matrix, if $AZ\geq kZ$ for a nonegative nonzero vector $Z$, then $k\leq\tau(A)$.
\end{lemma}

\begin{lemma}\label{lem:5.2}(\cite{HL}).
Let $x_{j}=(x_{j}(1),\cdots,x_{j}(n))^{T}\geq0$, $j\in\{1,2\cdots,m\}$, if $P_{j}>0$ and $\sum^{m}_{k=1}\frac{1}{P_{k}}\geq 1$, then we have
\begin{equation}\label{eq:5.1}
\sum\limits_{i=1}^{n}\prod_{j=1}^{m}x_{j}(i)\leq\prod_{j=1}^{m}\Big\{\sum\limits_{i=1}^{n}[x_{j}(i)]^{P_{j}}\Big\}^{\frac{1}{P_{j}}}.
\end{equation}
\end{lemma}

Next, according to these results, we expand the inequality
(\ref{eq:3.2}) of the Fan product of two matrices to the Fan product
of several matrices. One can obtain the following result:

\begin{theorem}\label{th:5.1}
For any matrices $A_{k}\in M_{n}$, and positive integers $P_{k}$
with $\sum^{m}_{k=1}\frac{1}{P_{k}}\geq 1$, $k\in\{1,2,\cdots,m\}$,
we have that
\begin{equation}\label{eq:5.2}
\tau(A_{1}\star A_{2}\cdots\star A_{m})\geq \min_{1\leq i\leq
n}\Big\{\prod^{m}_{k=1}A_{k}(i,i)-\prod^{m}_{k=1}[A_{k}(i,i)^{P_{k}}-\tau(A_{k}^{(P_{k})})]^\frac{1}{P_{k}}\Big\}.
\end{equation}
\end{theorem}

\textbf{Proof.} It is quite evident that the (\ref{eq:5.2}) holds with
the equality for $n=1$. Below we assume that $n\geq 2$.

\textbf{Case 1}. Let $A_{1}\star A_{2}\cdots\star A_{m}$ be an
irreducible nonsingular $M$-matrix, thus $A_{k}$ is irreducible,
$k\in\{1,2,\cdots,m\}$, we can obtain that $A_{k}^{(P_{k})}$ is also
irreducible. Let
$u_{k}^{(P_{k})}=(u_{k}(1)^{P_{k}},\cdots,u_{k}(n)^{P_{k}})^T>0$ be
a right Perron eigenvector of $A_{k}^{(P_{k})}$, and
$u_{k}=(u_{k}(1),\cdots,u_{k}(n))^T>0$, thus for any $i\in N$, we
have that
$$
A_{k}^{(P_{k})}u_{k}^{(P_{k})}=\tau (A_{k}^{(P_{k})})u_{k}^{(P_{k})},
$$
$$
A_{k}(i,i)^{P_{k}}u_{k}(i)^{P_{k}}-\sum\limits_{j\neq
i}|A_{k}(i,j)^{P_{k}}|u_{k}(j)^{P_{k}}=\tau
(A_{k}^{(P_{k})})u_{k}(i)^{(P_{k})},
$$
and
\begin{equation}\label{eq:5.3}
\sum\limits_{j\neq
i}|A_{k}(i,j)^{P_{k}}|u_{k}(j)^{P_{k}}=\Big(A_{k}(i,i)^{P_{k}}-\tau
(A_{k}^{(P_{k})})\Big)u_{k}(i)^{P_{k}}.
\end{equation}
Denote $C=A_{1}\star A_{2}\cdots\star A_{m}$, $Z=u_{1}\star
u_{2}\cdots\star u_{m}=(Z(1),\cdots,Z(n))^{T}>0$, thus
$Z(i)=\prod^{m}_{k=1}u_{k}(i)$. By the Lemma \ref{lem:5.2}
and (\ref{eq:5.3}), we get that
$$
\begin{array}{lll}
(CZ)_{i}=\Big(\prod^{m}_{k=1}A_{k}(i,i)\Big)Z(i)-\Big(\sum\limits_{j\neq
i}\prod^{m}_{k=1}|A_{k}(i,j)|\Big)Z(j)\\~~~~~~~~~
=\Big(\prod^{m}_{k=1}A_{k}(i,i)\Big)Z(i)-\sum\limits_{j\neq
i}\prod^{m}_{k=1}\Big(|A_{k}(i,j)|u_{k}(j)\Big)\\~~~~~~~~~
\geq\Big(\prod^{m}_{k=1}A_{k}(i,i)\Big)Z(i)-\prod^{m}_{k=1}\Big\{\sum\limits_{j\neq
i}[|A_{k}(i,j)|u_{k}(j)]^{(P_{k})}\Big\}^\frac{1}{P_{k}}~~\mathrm{(by~the~
equality~ (\ref{eq:5.3}))}\\~~~~~~~~~
=\Big(\prod^{m}_{k=1}A_{k}(i,i)\Big)Z(i)-\prod^{m}_{k=1}\Big\{[A_{k}(i,i)^{P_{k}}-\tau
(A_{k}^{(P_{k})})]u_{k}(i)^{P_{k}}\Big\}^\frac{1}{P_{k}}\\~~~~~~~~~
=\Big\{\prod^{m}_{k=1}A_{k}(i,i)-\prod^{m}_{k=1}[A_{k}(i,i)^{P_{k}}-\tau
(A_{k}^{(P_{k})})]\Big\}^\frac{1}{P_{k}}Z(i).
\end{array}
$$
According to the Lemma \ref{lem:5.1}, we obtain that
$$
\tau(A_{1}\star A_{2}\cdots\star A_{m})\geq \min_{1\leq i\leq
n}\Big\{\prod^{m}_{k=1}A_{k}(i,i)-\prod^{m}_{k=1}[A_{k}(i,i)^{P_{k}}-\tau(A_{k}^{(P_{k})})]^\frac{1}{P_{k}}\Big\}.
$$

\textbf{Case 2}. If $A_{1}\star A_{2}\cdots\star A_{m}$ is
reducible, where $A_{i}~(i=1,2,\cdots,m)$ are nonsingular
$M$-matrices. Similarly, let $P=(p_{ij})$ be the $n\times n$
permutation matrix with $p_{12}=p_{23}=\cdots=p_{n-1,n}=p_{n,1}=1$,
the remaining $p_{ij}$ zero, then $A_{k}-\varepsilon P$ is an
irreducible nonsingular $M$-matrix for any chosen positive real
number $\varepsilon$. Note that $A_{k}-\varepsilon P$ is a
continuous function on $\varepsilon$. Now we substitute
$A_{k}-\varepsilon P$ for $A_{k}$, in the previous Case 1, and then
letting $\varepsilon\rightarrow 0$, the result (\ref{eq:5.2})
follows by continuity. $\Box$

\textbf{Remark 4.1.}\label{re:5.1} If we take $m=2$ in Theorem
\ref{th:4.1}, one can obtain the following results:

\begin{itemize}
  \item If $p_{1}=p_{2}=1$, $A_{1}=A=(a_{ij})$, $A_{2}=B=(b_{ij})$,
we have that
$$
\tau(A\star B)\geq \min_{1\leq i\leq
n}\Big\{a_{ii}b_{ii}-(a_{ii}-\tau(A))(b_{ii}-\tau(B))\Big\},
$$
which is just the inequality (\ref{eq:3.2}).
  \item If $p_{1}=p_{2}=2$, $A_{1}=A=(a_{ij})$, $A_{2}=B=(b_{ij})$, then
\begin{equation}\label{eq:5.4}
\tau(A\star B)\geq \min_{1\leq i\leq
n}\Big\{a_{ii}b_{ii}-[a_{ii}^{2}-\tau(A\star
A)]^{\frac{1}{2}}[b_{ii}^{2}-\tau(B\star B)]^{\frac{1}{2}}\Big\}.
\end{equation}
In addition, by using the inequalities of arithmetic
and geometric means, we may obtain that
\[
a_{ii}^{2}\tau(B\star B)+b_{ii}^{2}\tau(A\star A)\geq
2a_{ii}b_{ii}[\tau(A\star A)\tau(B\star B)]^{\frac{1}{2}},
\]
so
\begin{equation}\label{eq:5.5}
(a_{ii}^{2}-\tau(A\star A))(b_{ii}^{2}-\tau(B\star B))\leq
\Big\{a_{ii}b_{ii}-[\tau(A\star A)\tau(B\star B)]^{\frac{1}{2}}\Big\}^{2}.
\end{equation}
Since for any $A,B\in M_n$, $\tau(A\star B)\geq \tau(A) \tau(B)$
(see \cite{QG} or (\ref{eq:3.1})), then, by (\ref{eq:5.5}), we have
that
 $$
a_{ii}b_{ii}-\Big[(a_{ii}^{2}-\tau(A\star A))(b_{ii}^{2}-\tau(B\star
B))\Big]^{\frac{1}{2}}\geq[\tau(A\star A)\tau(B\star
B)]^{\frac{1}{2}}\geq\tau(A)\tau(B).
$$
That is, the bound in (\ref{eq:5.2}) is better than the bound in
(\ref{eq:3.1}).
  \item If $p_{1}=1, p_{2}=2, A_{1}=A=(a_{ij}), A_{2}=B=(b_{ij})$,
  then we get
$$
\tau(A\star B)\geq \min_{1\leq i\leq
n}\Big\{a_{ii}b_{ii}-[a_{ii}-\tau(A)][b_{ii}^{2}-\tau(B\star
B)]^{\frac{1}{2}}\Big\}.
$$
\end{itemize}

\textbf{Acknowledgements}. \emph{The authors sincerely thank Prof.
Julio Moro and the reviewers and editor for their valuable and
detailed comments and suggestions on the manuscript of this paper,
which led to a substantial improvement on the presentation and
contents of this paper.}

\end{document}